\documentclass[letterpaper]{article}
\usepackage{aaai16}
\usepackage{times}
\usepackage{helvet}
\usepackage{courier}
\usepackage{url}
\usepackage{pdfsync}
\usepackage{CJKutf8}
\usepackage{xcolor}
\usepackage{subfigure}
\usepackage{listings}
\usepackage[cmex10]{amsmath}
\usepackage{amsthm}
\usepackage{graphicx}
\usepackage{multirow}
\usepackage{multicol}
\usepackage{bbding}
\usepackage{amsfonts,amssymb}
\usepackage{bm}
\usepackage{algorithm}
\usepackage{algorithmic}
\usepackage{epstopdf}
\usepackage{enumerate}
\usepackage[longnamesfirst]{natbib}
\newtheorem{theorem}{Theorem}

\newtheorem{proposition}{Proposition}

\newtheorem{Definition}{Definition}
\newtheorem{remark}{Remark}
\newcommand{\y}{\mathbf{y}}
\newcommand{\x}{\mathbf{x}}

\newcommand{\z}{\mathbf{z}}
\newcommand{\dd}{\mathbf{d}}

\renewcommand{\u}{\mathbf{u}}
\renewcommand{\v}{\mathbf{v}}

\renewcommand{\b}{\mathbf{b}}

\newcommand{\A}{\mathcal{A}}

\newcommand{\lbar}{\left\|}
\newcommand{\rbar}{\right\|}

\DeclareMathOperator*{\argmin}{argmin}
\begin{document}
%
\title{Relaxed Majorization-Minimization for Non-smooth and Non-convex Optimization}
\author{Chen Xu$^1$, Zhouchen Lin$^{1,2,}$\thanks{Corresponding author.}, Zhenyu Zhao$^3$, Hongbin Zha$^1$\\
	$^1$ Key Laboratory of Machine Perception (MOE), School of EECS, Peking University, P. R. China\\
	$^2$ Cooperative Medianet Innovation Center, Shanghai Jiaotong University, P. R. China\\
	$^3$Department of Mathematics, School of Science, National University of Defense Technology, P. R. China\\
	{\tt\small xuen@pku.edu.cn, zlin@pku.edu.cn, dwightzzy@gmail.com, zha@cis.pku.edu.cn}
}
\maketitle
\begin{abstract}
\begin{quote}
We propose a new majorization-minimization (MM) method for
non-smooth and non-convex programs, which is general enough to
include the existing MM methods. Besides the local majorization
condition, we only require that the difference between the
directional derivatives of the objective function and its
surrogate function vanishes \emph{when the number of iterations
approaches infinity}, which is a very weak condition. So our
method can use a surrogate function that directly approximates the
non-smooth objective function. In comparison, all the existing MM
methods construct the surrogate function by approximating the
smooth component of the objective function. We apply our relaxed
MM methods to the robust matrix factorization (RMF) problem with
different regularizations, where our locally majorant algorithm
shows great advantages over the state-of-the-art approaches for
RMF. This is the first algorithm for RMF ensuring, \emph{without
extra assumptions}, that any limit point of the iterates is a
stationary point. 
\end{quote}
\end{abstract}
\section{Introduction}
\vspace{-0.2em}
Consider the following  optimization problem:
\begin{equation}\label{nonsmoothnonconvex}
\vspace{-0.2em}
    \min_{\x \in \mathcal{C}} ~ f(\x),
    \vspace{-0.1em}
\end{equation}
where $\mathcal{C}$ is a closed convex subset in $\mathbb{R}^n$
and $f(\x): \mathbb{R}^n \rightarrow \mathbb{R}$ is a continuous
function bounded below, which could be non-smooth and non-convex.
Often, $f(x)$ can be split as:
\begin{equation}\label{splitting}
  \vspace{-0.2em}
     f(\x)= \tilde{f} (\x)+ \hat{f}(\x),
     \vspace{-0.1em}
\end{equation}
where $\tilde{f} (\x)$ is differentiable and $\hat{f}(\x)$ is
non-smooth\footnote{$\tilde{f}(\x)$ and $\hat{f}(\x)$ may
vanish.}. Such an optimization problem is ubiquitous, e.g., in statistics \citep{smoothingmethod},  computer vision and image processing \citep{nonsmoothproblem1,takeoalternating}, data mining and machine learning \citep{l1nfm,l21rnmf}.
There have been a variety of methods to tackle problem
\eqref{nonsmoothnonconvex}. Typical methods include subdifferential
\citep{nonsmooth}, bundle methods \citep{bundle}, gradient sampling \citep{gradientsampling}, 
smoothing methods \citep{smoothingmethod}, and
majorization-minimization (MM) \citep{mmtutorial}. In
this paper, we focus on the MM methods.
\vspace{-0.3em}
\subsection{Existing MM for Non-smooth and Non-convex Optimization}
\vspace{-0em}
\renewcommand{\algorithmicrequire}{\textbf{Input:}}
\renewcommand{\algorithmicensure}{\textbf{Output:}}
\begin{algorithm}[t]
\caption{Sketch of MM} \label{algorithmMM}
\begin{algorithmic}[1]
\REQUIRE$\x_0 \in \mathcal{C}$.
 \WHILE {not converged}
 \STATE Construct a surrogate function $g_k(\x)$ of $f(\x)$ at the current
iterate $\x_k$.
 \STATE Minimize the surrogate to get the next
iterate: $\x_{k+1}=\argmin_{\x \in \mathcal{C}}g_k(\x)$.
 \STATE $k \leftarrow k+1$.
 \ENDWHILE
\ENSURE The solution $\x_k$.
\end{algorithmic}
\end{algorithm}

\begin{table}[!t] \vspace{-0.5em}
\caption{Comparison of surrogate functions among existing MM
methods. $\mathbf{\#1}$ represents globally majorant MM in~\citep{mairalsurrogate}, $\mathbf{\#2}$ represents strongly convex MM in~\citep{mairalsurrogate}, $\mathbf{\#3}$ represents successive MM in~\citep{uppernonsmooth} and $\mathbf{\#4}$ represents our relaxed MM. In the second to fourth rows, $\times$ means that a function is not necessary to have the feature. In the last three rows, $\surd$ means that a function has the ability.}
   \renewcommand{\arraystretch}{1.3}
\label{summary}
\begin{center}
 \begin{tabular}{|c|c|c|c|c|}
\hline
 \footnotesize{Surrogate Functions}& \footnotesize{$\mathbf{\#1}$}& \footnotesize{$\mathbf{\#2}$ }& \footnotesize{$\mathbf{\#3}$} & \footnotesize{$\mathbf{\#4}$} 
    \\\hline\hline
\footnotesize{Globally Majorant} &$\surd$&$\times$&$\surd$&$\times$ \\ \hline
\footnotesize{Smoothness of Difference}  &\multirow{1}{*}{$\surd$}&\multirow{1}{*}{$\surd$}&\multirow{1}{*}{$\times$}&\multirow{1}{*}{$\times$} \\\hline
\footnotesize{Equality of Directional Derivative}  &\multirow{1}{*}{$\surd$}&\multirow{1}{*}{$\surd$}&\multirow{1}{*}{$\surd$}&\multirow{1}{*}{$\times$} \\\hline \hline
\multirow{1}{*}{\footnotesize{Approximate $\tilde{f}(\x)$} in \eqref{splitting}}  &\multirow{1}{*}{$\surd$}&\multirow{1}{*}{$\surd$}&\multirow{1}{*}{$\surd$}&\multirow{1}{*}{$\surd$} \\\hline
\footnotesize{Approximate $f(\x)$ in \eqref{nonsmoothnonconvex} ($\hat{f}\neq 0$)}  &\multirow{1}{*}{$\times$}&\multirow{1}{*}{$\times$}&\multirow{1}{*}{$\times$}&\multirow{1}{*}{$\surd$} \\ \hline
\multirow{1}{*}{\footnotesize{Sufficient Descent}}
&\multirow{1}{*}{$\times$}&\multirow{1}{*}{$\surd$}&\multirow{1}{*}{$\times$}&\multirow{1}{*}{$\surd$}\\ \hline
\end{tabular}
\end{center} \vspace{-3em}
\end{table}
\noindent MM has been successfully applied to a wide range of problems. \citet{mairalsurrogate} has given a comprehensive review on
MM. Conceptually, the MM methods consist of two steps (see
Algorithm \ref{algorithmMM}). First, construct a surrogate
function $g_k(\x)$ of $f(\x)$ at the current iterate $\x_k$.
Second, minimize the surrogate $g_k(\x)$ to update $\x$. The
choice of surrogate is critical for the efficiency of solving
\eqref{nonsmoothnonconvex} and also the quality of solution. The
most popular choice of surrogate is the class of ``first order
surrogates",  whose difference from the objective function  is
differentiable with a Lipschitz continuous gradient \citep{mairalsurrogate}. For non-smooth and non-convex objectives, to the best of our
knowledge, ``first order surrogates" are only used to approximate the
differentiable part $\tilde{f}(\x)$ of the objective $f(\x)$ in
\eqref{splitting}. More precisely, denoting $\tilde{g}_k
(\x)$ as an approximation of $\tilde{f} (\x)$ at iteration $k$,
the surrogate is
    \begin{equation} \label{smoothnonsmooth}
     g_k(\x)= \tilde{g}_k (\x)+ \hat{f}(\x).
   \end{equation}
Such a split approximation scheme has been successfully applied,
e.g., to minimizing the difference of convex functions
\citep{candesdcprogramming} and in the proximal splitting algorithm
\citep{nonconvexpg}. In parallel to \citep{mairalsurrogate}, \citet{uppernonsmooth} also showed that many
popular methods for minimizing non-smooth functions could be
regarded as MM methods. They proposed the block coordinate
descent method, where the traditional MM could be regarded as a
special case by gathering all variables in one block. Different
from \citep{mairalsurrogate}, they suggested using the directional
derivative to ensure the first order smoothness between the
objective and the surrogate, which is weaker than the condition in
\citep{mairalsurrogate} that the difference between the objective
and the surrogate should be smooth. However,
\citet{uppernonsmooth} only discussed the choice of surrogates by
approximating $f (\x)$ as \eqref{smoothnonsmooth}. \vspace{-0.2em}
\subsection{Contributions}
The contributions of this paper are as follows:
\begin{enumerate}[(a)]
 \item We further relax the condition on the
difference between the objective and the surrogate. We only
require that the directional derivative of the difference vanishes
\emph{when the number of iterations approaches infinity} (see
\eqref{surrogate3}). Our even weaker condition ensures that the
non-smooth and non-convex objective can be approximated directly.
Our relaxed MM is general enough to include the existing MM
methods \citep{mairalsurrogate,uppernonsmooth}.
 \item We also propose the
conditions ensuring that the iterates produced by our relaxed MM
converge to stationary points\footnote{For convenience, we say
that a sequence converges to stationary points meaning that any
limit point of the sequence is a stationary point.}, even for
general non-smooth and non-convex objectives.
\vspace{-0em}
\begin{figure}[t]
\centering
\includegraphics[width=0.7\linewidth]{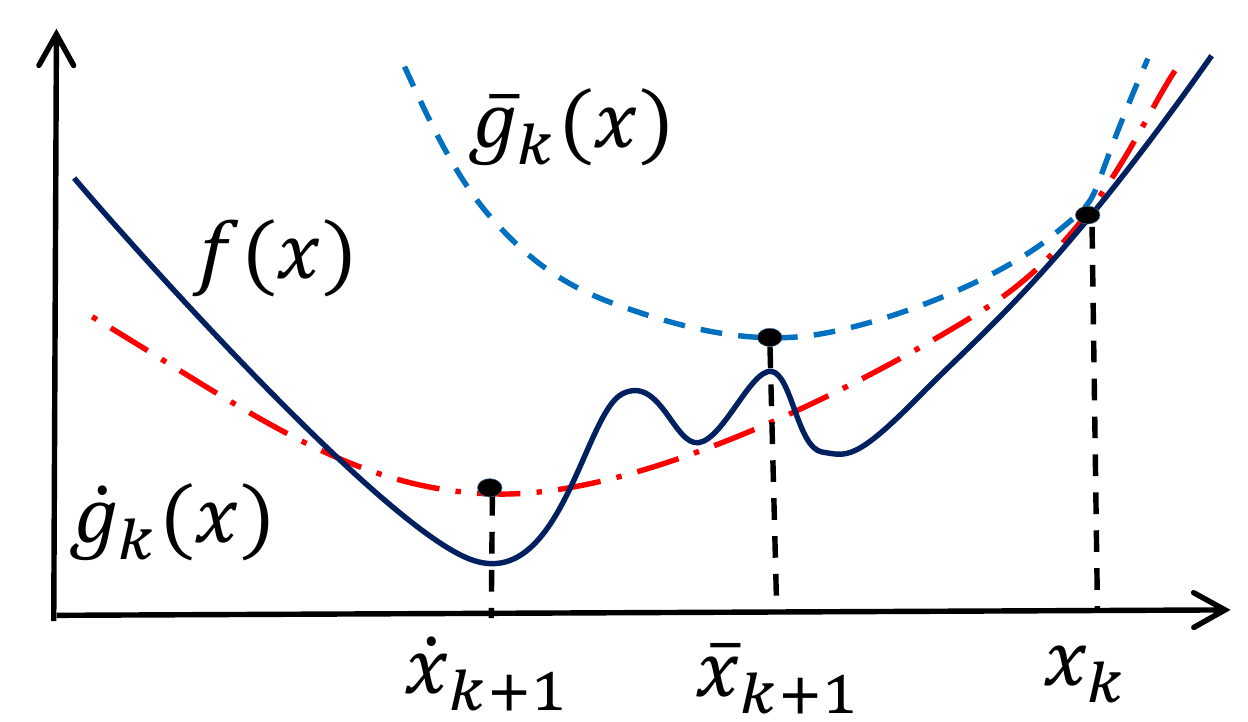}
\vspace{-1em}
\caption{Illustration of a locally majorant surrogate
$\dot{g}_k(\x)$ and a globally majorant surrogate $\bar{g}_k(\x)$.
Quite often, a globally majorant surrogate cannot approximate the
objective function well, thus giving worse solutions.}
\label{MM_illistration} \vspace{-1em}
\end{figure}  
 \item As a concrete example, we apply our relaxed MM to the robust matrix factorization
(RMF) problem with different regularizations. Experimental results
testify to the robustness and effectiveness of our locally majorant algorithm over the
state-of-the-art algorithms for RMF. To the best of our knowledge,
this is the first work that ensures convergence to stationary
points \emph{without extra assumptions}, as the objective and the
constructed surrogate naturally fulfills the convergence conditions
of our relaxed MM.
\vspace{-0.2em}
\end{enumerate}
Table \ref{summary} summarizes the differences between our relaxed
MM and the existing MM. We can see that ours is general enough to
include existing works \citep{mairalsurrogate,uppernonsmooth}. In
addition, it requires less smoothness on the difference between
the objective and the surrogate and has higher approximation
ability and better convergence property.

\section{Our Relaxed MM}\label{generalMM}
Before introducing our relaxed MM, we recall some definitions that
will be used later.
\begin{Definition}\label{def1}(\textbf{Sufficient Descent}) $\{f(\x_k)\}$ is said to have sufficient descent on the sequence $\{\x_k\}$ if there exists a constant $\alpha > 0$ such that:
      \begin{equation}
           f(\x_k) - f(\x_{k+1}) \geq \alpha \|\x_k - \x_{k+1}\|^2,\quad\forall k.
      \end{equation}
\end{Definition}
\begin{Definition}\label{def2}(\textbf{Directional Derivative} \citep[Chapter 6.1]{convexnonlinear}) The directional derivative of function $f(\x)$ in the feasible direction $\dd$ ($\x+\dd\in \mathcal{C}$) is defined as:
      \begin{equation}
           \nabla f(\x;\dd) = \liminf_{\theta \downarrow 0} \frac{f(\x+\theta \dd) - f(\x)}{ \theta}.
      \end{equation}
\end{Definition}

\begin{Definition}\label{def3} (\textbf{Stationary Point} \citep{uppernonsmooth}) A point $\x^*$ is a (minimizing) stationary  point of $f(\x)$ if $\nabla f(\x^*;\dd) \geq 0$ for all $\dd$ such that $\x^*+\dd \in \mathcal{C}$.
\end{Definition}
The surrogate function in our relaxed MM should satisfy the
following three conditions:
\begin{equation}\label{surrogate1}
  f(\x_k)=g_k(\x_k),
\end{equation}
\begin{equation}\label{surrogate2}
  f(\x_{k+1}) \leq  g_k(\x_{k+1}),\quad \mbox{(Locally Majorant)}
\end{equation}
\begin{equation} \label{surrogate3}
\begin{split}
 &\lim_{k \to \infty}  \left( \nabla f(\x_k; \dd)- \nabla g_k(\x_k;\dd) \right)= 0,  \\  & \forall \x_k+\dd \in
 \mathcal{C}. \quad \mbox{(Asymptotic Smoothness)}
 \end{split}
\end{equation}
By combining conditions \eqref{surrogate1} and \eqref{surrogate2}, we
have the non-increment property of MM:
\begin{equation}\label{descending}
     f(\x_{k+1}) \leq g_k(\x_{k+1})\leq g_k(\x_k)  = f(\x_k).
\end{equation}
However, we will show that with a careful choice of surrogate,
$\{f(\x_k)\}$ can have sufficient descent, which is stronger than
non-increment and is critical for proving the convergence of our
MM.

In the traditional MM, the global majorization condition
\citep{mairalsurrogate,uppernonsmooth} is assumed:
\begin{equation}\label{surrogate2s}
f(\x) \leq  g_k(\x),  \quad \forall \x \in \mathcal{C},\quad
\mbox{(Globally Majorant)}
\end{equation}
which also results in the non-increment of the objective function,
i.e., \eqref{descending}. However, a globally majorant surrogate
cannot approximate the object well (Fig.~\ref{MM_illistration}).
Moreover, the step length between successive iterates may
be too small. So a globally majorant surrogate is likely to
produce an inferior solution and converges slower than a locally
majorant one (\citep{mairalsurrogate} and our experiments).

Condition \eqref{surrogate3} requires very weak first order
smoothness of the difference between $g_k(\x)$ and $f(\x)$. It is
weaker than that in \cite[Assumption~1.(A3)]{uppernonsmooth},
which requires that the directional derivatives of $g_k(\x)$ and
$f(\x)$ are equal at every $\x_k$. Here we only require the
equality when the number of iterations goes to infinity, which provides more flexibility in constructing the surrogate function.
When $f(\cdot)$ and $g_k(\cdot)$ are both smooth, condition
\eqref{surrogate3} can be satisfied when the two hold the same gradient at
$\x_k$. For non-smooth functions, condition
\eqref{surrogate3} can be enforced on the differentiable part
$\tilde{f}(\x)$ of $f(\x)$. These two cases have been discussed in the literature \citep{mairalsurrogate,uppernonsmooth}. If $\tilde{f}(\x)$ vanishes, the case not yet discussed in the literature, the
condition can still be fulfilled by approximating the whole
objective function $f(\x)$ directly, as long as the resulted
surrogate satisfies certain properties, as stated
below\footnote{The proofs of the results in this paper can be
found in Supplementary Materials.}.

\begin{proposition}\label{proposition1}
Assume that $\exists~ K>0,~+ \infty
> \gamma_u$, $\gamma_l
>0$, and $\epsilon >0$, such
that
\begin{equation}\label{eqn:upper-and-lower-bound}
\hat{g}_k(\x)+\gamma_u \|\x-\x_k\|^2_2 \geq f(\x) \geq
\hat{g}_k(\x)- \gamma_l \|\x-\x_k\|^2_2
\end{equation}
holds for all $k\geq K$ and $\x \in \mathcal{C}$ such that
$\|\x-\x_k\|\leq \epsilon$, where the equality holds if and only
if $\x=\x_k$. Then condition \eqref{surrogate3} holds for
$g_k(\x)= \hat{g}_k(\x)+\gamma_u \|\x-\x_k\|^2_2$.
\end{proposition}
We make two remarks. First, for sufficiently large $\gamma_u$ and
$\gamma_l$ the inequality \eqref{eqn:upper-and-lower-bound}
naturally holds. However, a larger $\gamma_u$ leads to slower
convergence. Fortunately, as we only require the inequality to
hold for sufficiently large $k$, adaptively increasing $\gamma_u$
from a small value is allowed and also beneficial. In many cases,
the bounds of $\gamma_u$ and $\gamma_l$ may be deduced from the
objective function. Second, the above proposition does not specify
any smoothness property on either $\hat{g}_k(\x)$ or the
difference $g_k(\x)- f(\x)$. It is not odd to add the proximal
term $\|\x-\x_k\|^2_2$, which has been widely used, e.g.
in proximal splitting algorithm
\citep{nonconvexpg} and alternating direction method of multipliers (ADMM) \citep{ladmpsap}.

For general non-convex and non-smooth problems, proving the
convergence to a global (or local) minimum is out of reach,
classical analysis focuses on converging to stationary points
instead. For general MM methods, since only non-increment property
is ensured, even the convergence to stationary points cannot be
guaranteed. Mairal~\citep{mairalsurrogate} proposed using strongly
convex surrogates and thus proved that the iterates of
corresponding MM converge to stationary points. Here we have
similar results, as stated below.

\begin{theorem}\label{theorem} \textbf{(Convergence)}
Assume that the surrogate $g_k(\x)$ satisfies \eqref{surrogate1}
and \eqref{surrogate2}, and further is strongly convex,
then the sequence $\{f(\x_k)\}$ has sufficient
descent. If $g_k(\x)$ further satisfies \eqref{surrogate3} and
$\{\x_k\}$ is bounded, then the sequence
$\{\x_k\}$ converges to stationary points.
\end{theorem}
\begin{remark}If $g_k(\x)=\dot{g}_k(\x)+\rho/2 \|\x-\x_k\|_2^2$, where
$\dot{g}_k(\x)$ is locally majorant (not necessarily convex) as
\eqref{surrogate2} and $\rho > 0$, then the strongly convex
condition can be removed and the same convergence result holds.
\end{remark}

In the next section, we  will give a concrete example on how to
construct appropriate surrogates for the RMF problem.

\section{Solving Robust Matrix Factorization by Relaxed MM} \label{constructsurrogate}
Matrix factorization is to factorize a matrix, which usually has
missing values and noises, into two matrices. It is widely used for structure from motion \citep{rsfm}, clustering \citep{l21rnmf}, dictionary learning \citep{onlinemarial}, etc. Normally,
people aim at minimizing the error between the given matrix and
the product of two matrices at the observed entries, measured in
squared $\ell_2$ norm. Such models are fragile to outliers.
Recently, $\ell_1$-norm has been suggested to measure the error for
enhancing robustness. Such models are thus called robust matrix
factorization (RMF). Their formulation is as follows:
  \begin{equation} \label{eq:original}
    \min _{U\in \mathcal{C}_u,V\in \mathcal{C}_v} \|W\odot(M-UV^T)\|_1+ R_u(U)+R_v(V),
  \end{equation}
where $M\in \mathbb{R}^{m \times n}$ is the observed matrix and
$\|\cdot\|_1$ is the $\ell_1$-norm, namely the sum of absolute
values in a matrix. $U\in \mathbb{R}^{m \times r}$ and $V \in
\mathbb{R}^{n \times r}$  are the unknown factor matrices. $W$ is
the 0-1 binary mask with the same size as $M$. The entry value $0$
means that the corresponding entry in $M$ is missing, and $1$
otherwise. The operator $\odot$ is the Hadamard entry-wise
product. $\mathcal{C}_u \subseteq \mathbb{R}^{m \times r}$ and
$\mathcal{C}_v \subseteq \mathbb{R}^{n \times r}$ are some closed
convex sets, e.g., non-negative cones or balls in some norm.
$R_u(U)$ and $R_v(V)$ represent some convex regularizations,
e.g., $\ell_1$-norm, squared Frobenius norm, or elastic net. By
combining different constraints and regularizations, we can
 get variants of RMF, e.g., low-rank matrix recovery
\citep{takeoalternating}, non-negative matrix factorization (NMF)
\citep{nmf}, and dictionary learning \citep{onlinemarial}.

Suppose that we have obtained $(U_k,V_k)$ at the $k$-th iteration.
We split $(U,V)$ as the sum of $(U_k,V_k)$ and the unknown
increment $(\Delta U, \Delta V)$:
    \begin{equation} \label{eq:update}
          (U,V)=(U_k,V_k)+(\Delta U,\Delta V).
 \end{equation}
Then \eqref{eq:original} can be rewritten as:
 \begin{equation} \label{eq:linear}
 \begin{split}
 &\min_{\Delta U + U_k \in \mathcal{C}_u,  \Delta V+ V_k \in \mathcal{C}_v} F_k(\Delta U,\Delta V)=\|W\odot(M-(U_k+ \\&\Delta U)(V_k^T+\Delta V)^T)\|_1  +R_u(U_k+\Delta U) +R_v(V_k+\Delta V).
 \end{split}
 \end{equation}

Now we aim at finding an increment $(\Delta U,\Delta V)$ such that
the objective function decreases properly. However, problem
\eqref{eq:linear} is not easier than the original problem
\eqref{eq:original}. Inspired by MM, we try to approximate
\eqref{eq:linear} with a convex surrogate. By the triangular
inequality of norms, we have the following inequality:
 \begin{equation} \label{eq:in1}
  \begin{split}
   & F_k(\Delta U,\Delta V) \leq \|W\odot(M-U_kV_k^T-\Delta UV_k^T-U_k\Delta V^T)\|_1 \\ & +\|W \odot (\Delta U \Delta V ^T)\|_1 + R_u(U_k+\Delta U) +R_v(V_k+\Delta V),
  \end{split}
\end{equation}
where the term $\|W \odot (\Delta U \Delta V ^T)\|_1$ can be further
 approximated by $\rho_u/2\|U\|^2_F+\rho_v/2\|V\|^2_F$, in which
$\rho_u$ and $\rho_v$ are some positive constants. Denoting
\begin{equation}
\begin{split}
&\hat{G}_k(\Delta U,\Delta V)=\|W\odot(M-U_kY_k^T-\Delta
UV_k^T \\&-U_k\Delta V^T)\|_1 +R_u(U_k+\Delta U)+R_v(V_k+\Delta V),
\end{split}
\end{equation}
we have a surrogate function of $F_k(\Delta U, \Delta V)$ as
follows:
   \begin{equation} \label{eq:surrogate}
   \begin{split}
     G_k(\Delta U, \Delta V)& = \hat{G}_k(\Delta U, \Delta V) + \frac{\rho_u}{2}\|\Delta U\|_F^2 + \frac{\rho_v}{2}\|\Delta V\|_F^2,  \\ & \mbox{s.t. } \Delta U + U_k \in \mathcal{C}_u,  \Delta V+ V_k \in
     \mathcal{C}_v.
     \end{split}
   \end{equation}

Denoting $\# W_{(i,.)}$ and $\# W_{(.,j)}$ as the number of
observed entries in the corresponding column and row of $M$,
respectively, and $\epsilon
>0 $ as any positive scalar, we have the following proposition.

\begin{proposition} \label{proposition2}
  $\hat{G}_k(\Delta U, \Delta V) + \bar{\rho}_u/2\|\Delta U\|_F^2 + \bar{\rho}_v/2\|\Delta V\|_F^2 \geq F_k(\Delta U, \Delta V) \geq \hat{G}_k(\Delta U, \Delta V) - \bar{\rho}_u/2\|\Delta U\|_F^2 - \bar{\rho}_v/2\|\Delta V\|_F^2$ holds for all possible $(\Delta U, \Delta V)$ and the equality holds if and only if $(\Delta U, \Delta V)=\mathbf{(0, 0)}$, where $\bar{\rho}_u=\max \{\# W_{(i,.)},i=1,\ldots,m\}+\epsilon$, and $\bar{\rho}_v=\max \{\# W_{(.,j)},j=1,\ldots,n\}+\epsilon$.
\end{proposition}
By choosing $\rho_u$ and $\rho_v$ in different ways, we have two
versions of relaxed MM for RMF: RMF by globally majorant MM
(RMF-GMMM for short) and RMF by locally majorant MM (RMF-LMMM for
short). In RMF-GMMM, $\rho_u$ and $\rho_v$ are fixed to be
$\bar{\rho}_u$ and $\bar{\rho}_v$ in Proposition
\ref{proposition2}, respectively, throughout the iterations. In
RMF-LMMM, $\rho_u$ and $\rho_v$ are instead initialized with
relatively small values and then increase gradually, using the
line search technique in \citep{APG} to ensure the
locally majorant condition \eqref{surrogate2}. $\rho_u$ and
$\rho_v$ eventually reach the upper bounds $\bar{\rho}_u$ and
$\bar{\rho}_v$ in Proposition \ref{proposition2}, respectively. As
we will show, RMF-LMMM significantly outperforms RMF-GMMM in all
our experiments, in both convergence speed and quality of
solution.

Since the chosen surrogate $G_k$ naturally fulfills the conditions
in Theorem~\ref{theorem}, we have the following convergence result
for RMF solved by relax MM.

\begin{theorem}
 By minimizing \eqref{eq:surrogate} and updating $(U,V)$ according to \eqref{eq:update}, the sequence $\{F(U_k,V_k)\}$ has sufficient descent and the sequence $\{(U_k,V_k)\}$ converges to stationary points.
\end{theorem}
To the best of our knowledge, this is the first convergence
guarantee for variants of RMF without extra assumptions. In the
following, we will give two examples of RMF in
\eqref{eq:original}.

\subsection{Two Variants of RMF}
\vspace{0.2em}
 \textbf{Low Rank Matrix Recovery} exploits the
fact $r \ll \min(m,n)$ to recover the intrinsic low rank data
from the measurement matrix with missing data. When the error is
measured by the squared Frobenius norm, many algorithms have been
proposed \citep{dampednewton,largescale}. For robustness,  \citet{takeoalternating} proposed to adopt the
$\ell_1$-norm. They minimized $U$ and $V$ alternatively, which
could easily get stuck at non-stationary points
\citep{Bertsekas:nonlinear}. So \citet{l1wiberg} represented $V$
implicitly with $U$ and extended the Wiberg Algorithm to
$\ell_1$-norm. 
They only proved the convergence of the objective function value, not the sequence $\{(U_k,V_k)\}$ itself. Moreover, they had to assume that the dependence of $V$ on $U$ is differentiable, which is unlikely to hold everywhere. Additionally, as it unfolds matrix $U$ into a vector and adopt an () its memory requirement is very high, which prevents it from large scale computation.
Recently, ADMM was used for matrix recovery. By assuming
that the variables are bounded and convergent, \citet{lmafit} proved that any accumulation point of their
algorithm is the Karush-Kuhn-Tucker (KKT) point. However, the
method was only able to handle outliers with magnitudes comparable
to the low rank matrix \citep{lmafit}. Moreover, the penalty parameter was fixed, which was not easy to tune for fast convergence. Some researches \citep{regl1,unifying}
further extended that by adding different regularizations on
$(U,V)$ and achieved state-of-the-art performance. As the
convergence analysis in \citep{lmafit} cannot be directly extended,
it remains unknown whether the iterates converge to KKT points. In
this paper, we adopt the same formulation as \citep{unifying}:
\begin{equation}\label{lowrankmatrix}
        \vspace{-0.3em}
    \min _{U,V} \|W\odot(M-UV^T)\|_1+
    \frac{\lambda_u}{2}\|U\|_F^2+\frac{\lambda_v}{2}\|V\|_F^2,
            \vspace{-0em}
\end{equation}
where the regularizers $\|U\|_F^2$ and $\|V\|_F^2$ are for
reducing the solution space \citep{dampednewton}. \vspace{0.4em}

\noindent\textbf{Non-negative Matrix Factorization} (NMF) has been popular
since the seminal work of \citet{nmf}. 
\citet{l21rnmf} extended the squared $\ell_2$ norm to the
$\ell_{21}$-norm for robustness. Recently, \citet{l1nfm}
further introduced the $\ell_1$-norm to handle outliers in
non-negative dictionary learning, resulting in the following
model:
\begin{equation}\label{eqn:NMF}
    \min _{U \geq 0,V \geq 0} \|M-UV^T\|_1+ \frac{\lambda_u}{2}\|U\|_F^2+\lambda_v\|V\|_1,
\end{equation}
where $\|U\|_F^2$ is added in to avoid the trivial solution and
$\|V\|_1$ is to induce sparsity. All the three NMF models use
multiplicative updating schemes, which only differ in the weights
used. The multiplicative updating scheme is intrinsically a
globally majorant MM. Assuming that the iterates converge, they
proved that the limit of sequence is a stationary point. However,
 \citet{accnmf} pointed out that with such a
multiplicative updating scheme is hard to reach the convergence
condition even on toy data.
\subsection{Minimizing the Surrogate Function} \label{minimize the surrogate}
Now we show how to find the minimizer of the convex $G_k(\Delta U,
\Delta V)$ in \eqref{eq:surrogate}. This can be easily done by
using the linearized alternating direction method with parallel
splitting and adaptive penalty (LADMPSAP) \citep{ladmpsap}.
LADMPSAP fits for solving the following linearly constrained
separable convex programs:
\begin{equation}
\min \limits_{\x_1,\cdots,\x_n} \sum\limits_{j=1}^n f_j(\x_j),\quad
s.t.\quad
\sum \limits_{j=1}^n \mathcal{A}_j(\x_j)=\mathbf{b},\label{eq:separable constraint}
        \vspace{-0em}
\end{equation}
where $\x_j$ and $\mathbf{b}$ could be either vectors or matrices,
$f_j$ is a proper convex function, and $\mathcal{A}_j$ is a linear
mapping. To apply LADMPSAP, we first introduce an auxiliary matrix
$E$ such that $E =M-U_kY_k^T-\Delta U V_k^T-U_k\Delta V^T$. Then
minimizing $G_k(\Delta U, \Delta V)$ in \eqref{eq:surrogate} can
be transformed into:
  \begin{equation}\label{eq:eqconstraint}
   \begin{aligned}
& \min_{E,\Delta U,\Delta V}  \|W \odot E\|_1
 \\&+\left(\frac{\rho_u}{2}\|\Delta U\|_F^2+R_u(U_k+\Delta U)+ \delta_{\mathcal{C}_u}(U_k+\Delta U)\right)
 \\&+ \left(\frac{\rho_v}{2}\|\Delta V\|_F^2+R_v(V_k+\Delta V)+ \delta_{\mathcal{C}_v}(V_k+\Delta V)\right), \\
                 & \text{s.t.}\quad E+\Delta U V_k^T+U_k\Delta V^T=M-U_kY_k^T,
          \end{aligned}
              \end{equation}
where the indicator function $\delta_\mathcal{C}(\x): \mathbb{R}^p
\rightarrow \mathbb{R}$ is defined as:
  \begin{equation}
     \delta_{\mathcal{C}}(\x)=\left\{ \begin{array}{ll}
     0,  &\quad  \mbox{if } \x \in \mathcal{C},  \\
     +\infty, &\quad \mbox{otherwise.}
     \end{array}
        \right.
\end{equation}
Then problem \eqref{eq:eqconstraint} naturally fits into the model
problem \eqref{eq:separable constraint}. For more details, please
refer to Supplementary Materials.

\section{Experiments} \label{experiment}
In this section, we compare our relaxed MM algorithms with state-of-the-art RMF algorithms: UNuBi \citep{unifying} for low rank matrix recovery and $\ell_1$-NMF \citep{l1nfm} for robust NMF. The code of UNuBi \citep{unifying} was kindly provided by the authors. We implemented the code of $\ell_1$-NMF \citep{l1nfm} ourselves.
\subsection{Synthetic Data}
We first conduct experiments on synthetic data. Here we set the
regularization parameters $\lambda_u=\lambda_v=20/(m+n)$ and stop our relaxed MM algorithms
when the relative change in the objective function is less than
$10^{-4}$.
\begin{figure}[tb]
\centering
\begin{tabular}{ccc}
{\hspace{-13pt}} \includegraphics[width=0.48\linewidth,height=0.45\linewidth]{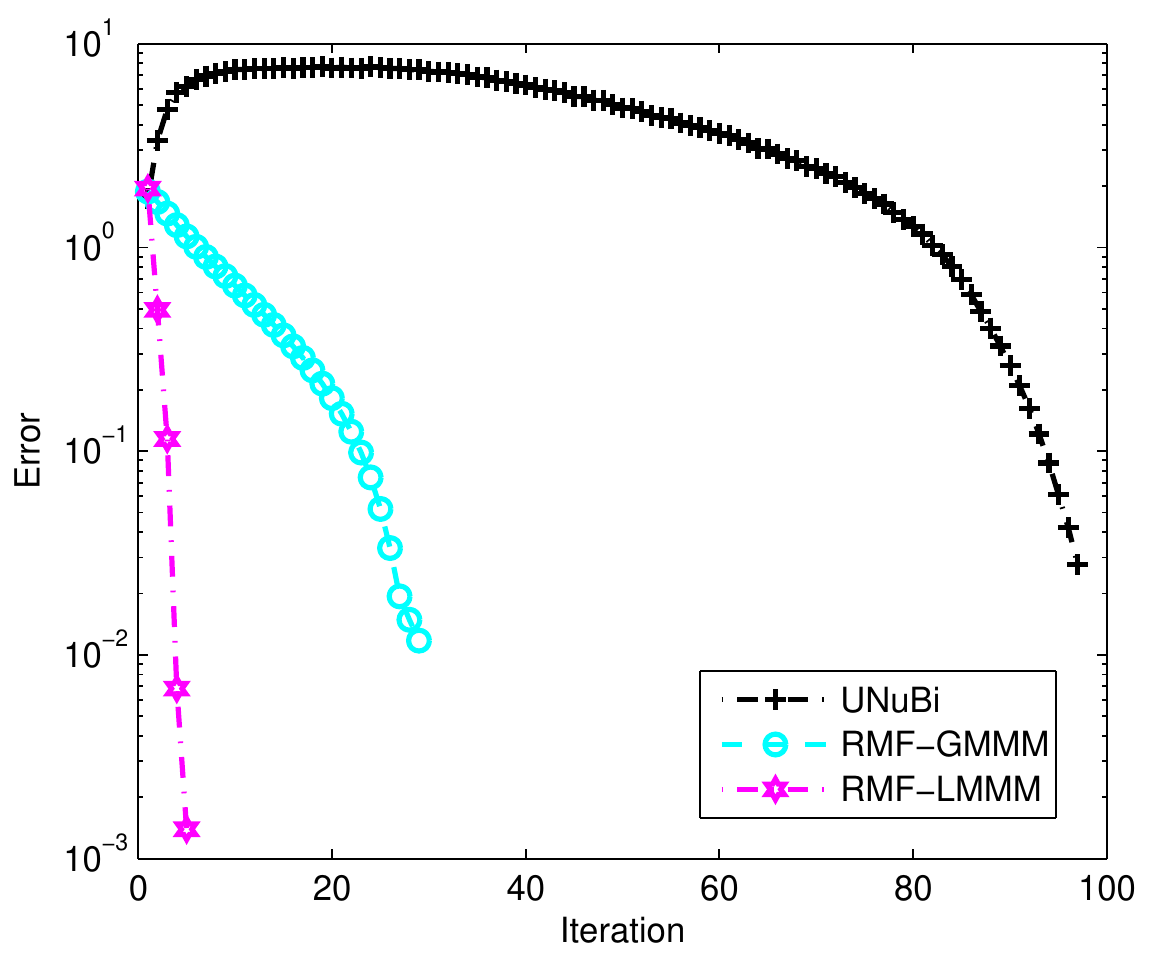}&~&
{\hspace{-13pt}} \includegraphics[width=0.48\linewidth,,height=0.45\linewidth]{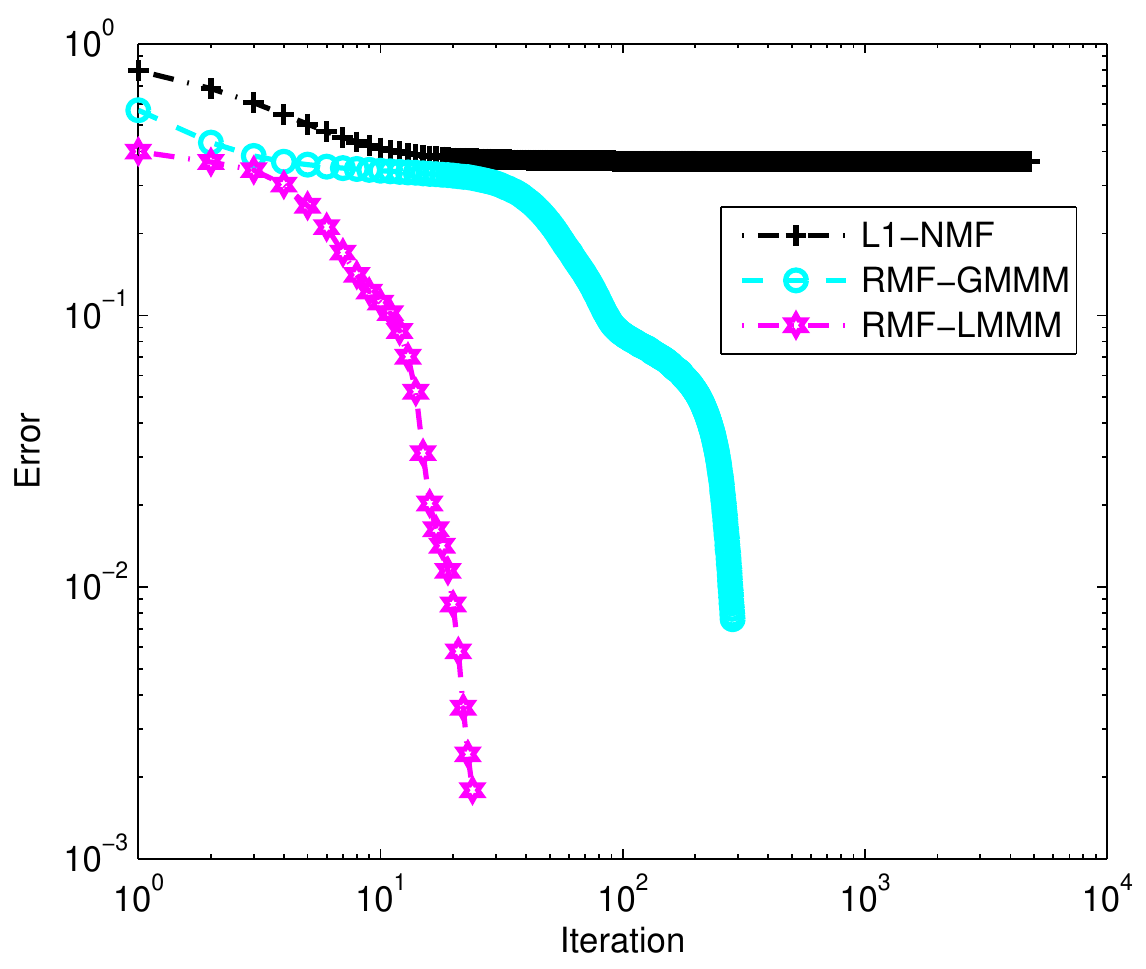}
{\vspace{-3.5pt}}\\
  \scriptsize{\hspace{-13pt}} (a) Low-rank Matrix Recovery &  &
  \scriptsize{\hspace{-13pt}} (b) Non-negative Matrix Factorization\\
  {\vspace{-1.7em}}
\end{tabular}
\caption{Iteration number versus relative error in log-$10$ scale
on synthetic data. (a) The locally majorant MM, RMF-LMMM, gets
better solution in much less iterations than the globally
majorant, RMF-GMMM, and the state-of-the-art algorithm, UNuBi
\citep{unifying}. (b) RMF-LMMM outperforms the globally
majorant MM, RMF-GMMM, and the state-of-the-art algorithm (also globally
majorant),
$\ell_1$-NMF \citep{l1nfm}. The iteration number in (b) is in
log-$10$ scale. }
\label{synthetic}
\end{figure}

\vspace{0.4em}
\noindent \textbf{Low Rank Matrix Recovery}: We generate a data matrix
$M=U_0V_0^T$, where $U_0 \in \mathbb{R}^{500\times 10}$ and $V_0
\in \mathbb{R}^{500\times 10}$ are sampled i.i.d. from a Gaussian
distribution $\mathcal{N}(0,1)$. We additionally corrupt $40\%$
entries of $M$ with outliers uniformly distributed in $[-10,10]$
and choose $W$ with $80\%$ data missing. The positions of both
outliers and missing data are chosen uniformly at random.  We initialize
all compared algorithms with the rank-$r$
truncation of the singular-value decomposition of $W\odot M$. The
performance is evaluated by measuring relative error with the
ground truth: $\|U_{est}V_{est}^T-U_0V_0^T\|_1/(mn)$, where
$U_{est}$ and $V_{est}$ are the estimated matrices. The results
are shown in Fig. \ref{synthetic}(a), where RMF-LMMM reaches the lowest
relative error in much less iterations.

\noindent\textbf{Non-negative Matrix Factorization}: We generate a data
matrix $M=U_0V_0^T$, where $U_0 \in \mathbb{R}^{500\times 10}$ and
$V_0 \in \mathbb{R}^{500\times 10}$ are sampled i.i.d. from a
uniform distribution $\mathcal{U}(0,1)$. For sparsity, we further
randomly set $30\%$ entries of $V$ as $0$. We further corrupt
$40\%$ entries of $M$ with outliers uniformly distributed in
$[0,10]$.  All the algorithms are initialized with the same
non-negative random matrix. The results are shown in Fig.
\ref{synthetic}(b), where RMF-LMMM also gets the best result with
much less iterations. $\ell_1$-NMF tends to be stagnant and cannot
approach a high precision solution even after $5000$ iterations.
\subsection{Real Data}
In this subsection, we conduct experiment on real data. Since there is no ground truth, we measure the relative error by $ \|W \odot (M_{est}-M)\|_1/ \# W$, where $\# W$ is the number of observed entries. For NMF, $W$ becomes an all-one matrix.

\begin{figure}[tb]
\begin{center}
\begin{tabular}{cccc}
{\hspace{-13pt}} \includegraphics[width=0.25\linewidth, height=0.25\linewidth]{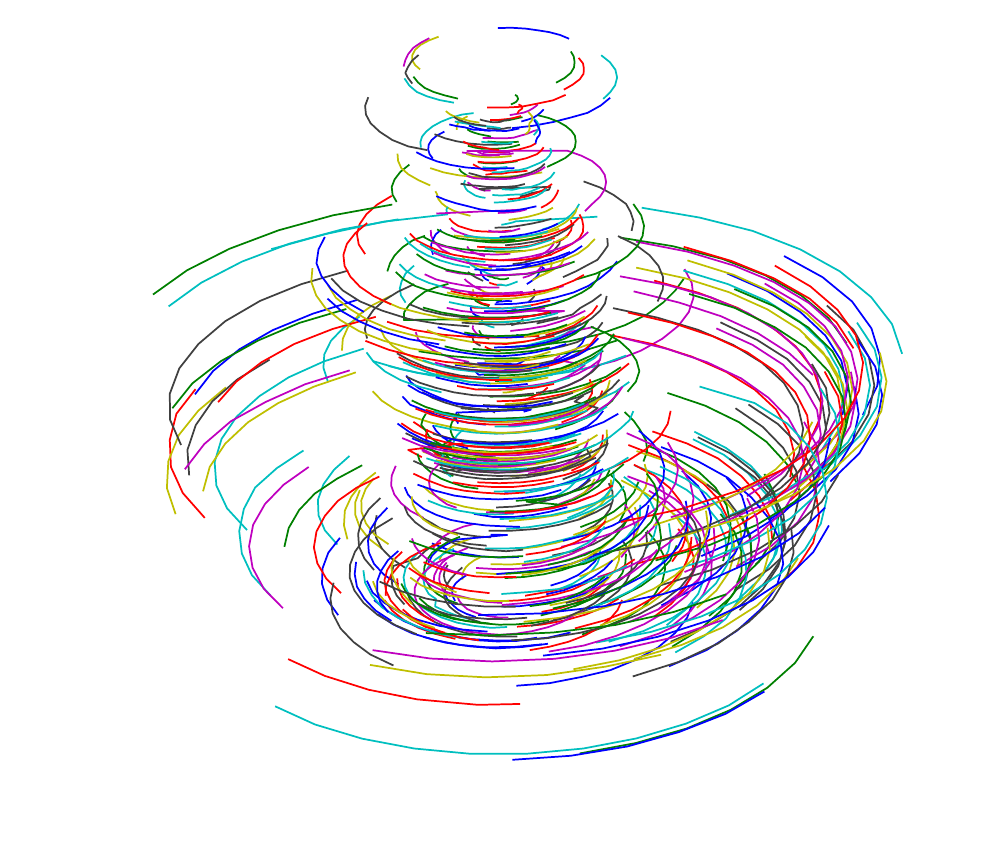}&
{\hspace{-13pt}} \includegraphics[width=0.25\linewidth,height=0.25\linewidth]{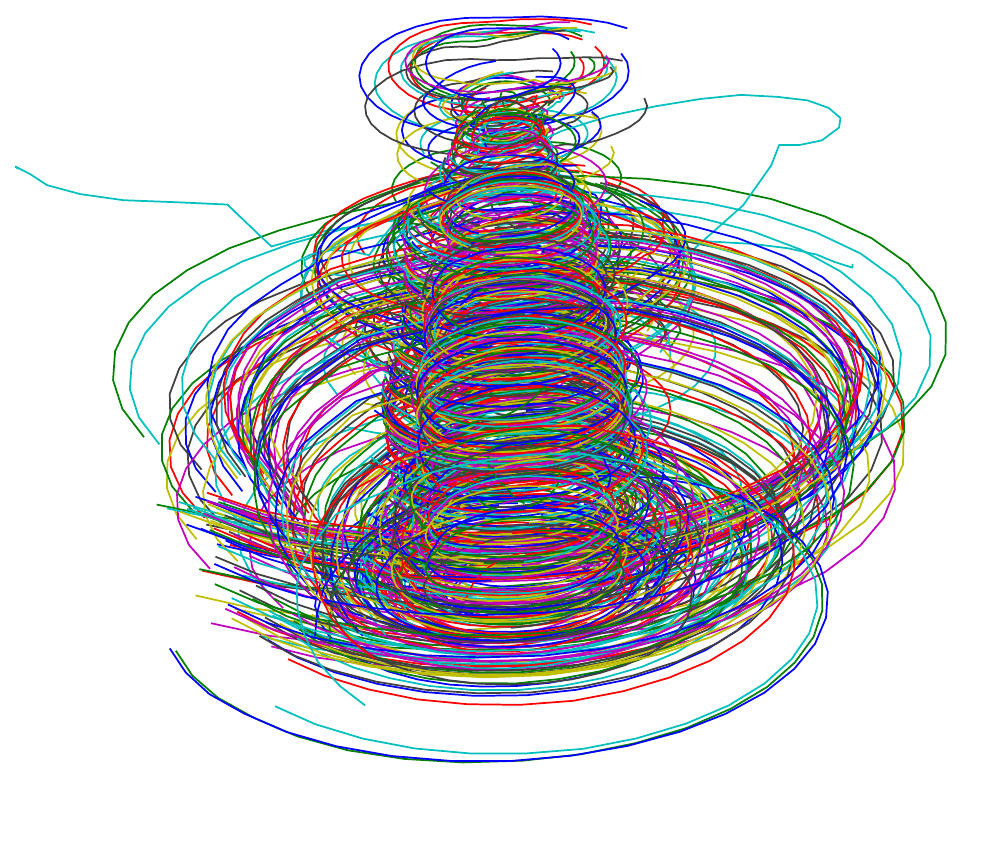}&
{\hspace{-13pt}} \includegraphics[width=0.25\linewidth,height=0.25\linewidth]{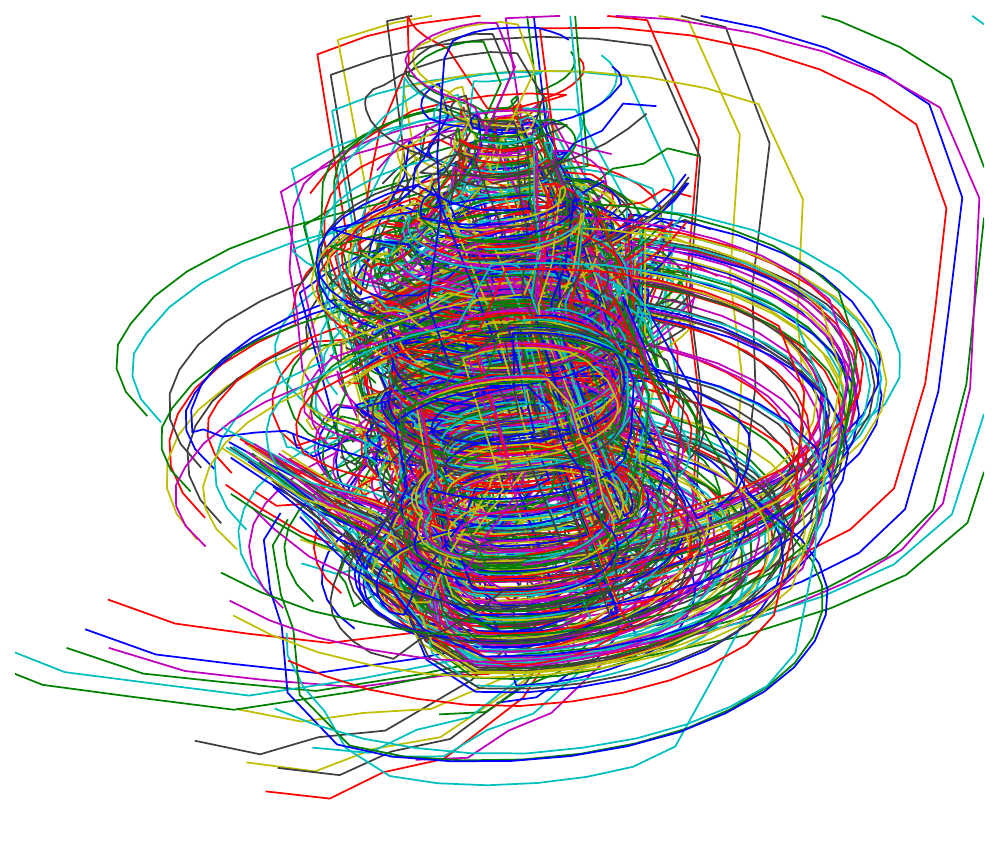}&
{\hspace{-13pt}} \includegraphics[width=0.25\linewidth, height=0.25\linewidth]{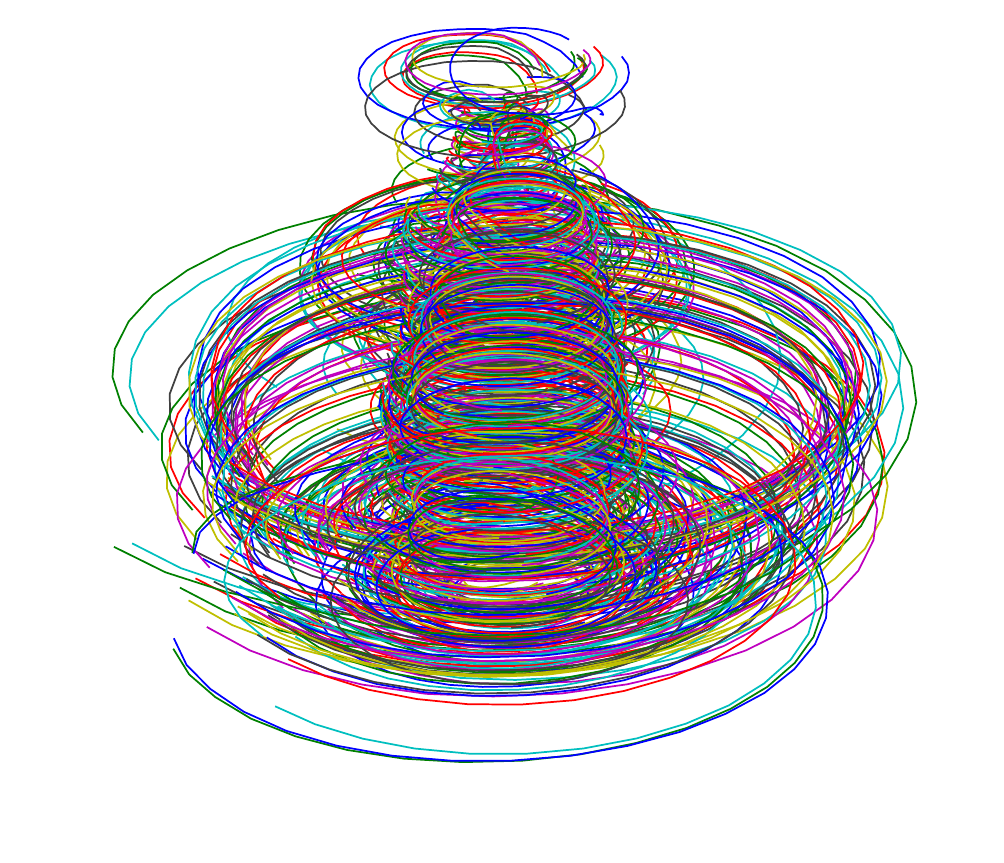}
{\vspace{-5pt}}\\
  \scriptsize{\hspace{-13pt}} (a) Raw Data &
  \scriptsize {\hspace{-13pt}} (b) UNuBi &
  \scriptsize {\hspace{-13pt}} (c) RMF-GMMM &
  \scriptsize {\hspace{-13pt}} (d) RMF-LMMM
{\vspace{-4pt}}\\
\scriptsize{\hspace{-13pt}}      &
  \scriptsize{\hspace{-13pt}}   Error=$0.336$ &
    \scriptsize{\hspace{-13pt}}  Error=$0.488$ &
  \scriptsize {\hspace{-13pt}}  Error=$\mathbf{0.322}$
\end{tabular}
{\vspace{-10pt}}\\
\end{center}
\caption{Original incomplete and recovered data of the Dinosaur
sequence. (a) Raw input tracks. (b-d) Full tracks reconstructed by
UNuBi~\citep{unifying}, RMF-GMMM, and RMF-LMMM, respectively.
 The relative error is presented below the tracks.} \vspace{-1em}
 \label{affine}
\end{figure}

\vspace{0.4em}
\noindent  \textbf{Low Rank Matrix Recovery}: \citet{rsfm} first modelled the affine rigid structure
from motion as a rank-$4$ matrix recovery problem. Here we use the
famous Oxford Dinosaur sequence  \footnote{\url{
http://www.robots.ox.ac.uk/~vgg/data1.html}},
which consists of $36$ images with a resolution of $720 \times
576$ pixels. We pick out a portion of the raw feature points which
are observed by at least $6$ views (Fig. \ref{affine}(a)). The
observed matrix is of size $72 \times 557$ with a missing data
ratio $79.5\%$ and shows a band diagonal pattern. We register the
image origin to the image center, $(360,288)$. We adopt
 the same initialization and parameter setting as the synthetic data above.

Figures \ref{affine}(b)-(d) show the full tracks reconstructed by
all algorithms. As the dinosaur sequence is taken on a turntable,
all the tracks are supposed to be circular. Among them, the tracks
reconstructed by RMF-GMMM are the most inferior. UNuBi gives
reasonably good results. However, most of the reconstructed tracks
in large radii do not appear closed. Some tracks in the upper part
are not reconstructed well either, including one obvious failure.
In contrast, almost all the tracks reconstructed by RMF-LMMM are
circular and appear closed, which are the most visually plausible.
The lowest relative error also confirms the effectiveness of
RMF-LMMM.

\vspace{0.4em}
\noindent \textbf{Non-negative Matrix Factorization}: We test the
performance of robust NMF by clustering \citep{l21rnmf,l1nfm}. The
experiments are conducted on four benchmark datasets of face
images, which includes: AT\&T,
UMIST, a subset of PIE
\footnote{\url{http://www.zjucadcg.cn/dengcai/Data/data.html}} and
a subset of AR
\footnote{\url{http://www2.ece.ohio-state.edu/~aleix/ARdatabase.html}}.
We use the first $10$ images in each class for PIE and the first
$13$ images for AR. The descriptions of the datasets are summarized
in the second row of Table \ref{Clustering results}. The
evaluation metrics we use here are accuracy (ACC), normalized
mutual information (NMI) and purity (PUR) \citep{l21rnmf,l1nfm}.
We change the regularization parameter $\lambda_u$ to $2000/(m+n)$
and maintain $\lambda_v$ as $20/(m+n)$. The number $r$ of clusters
is equal to the number of classes in each dataset. We adopt the
initializations in \citep{l21rnmf}. Firstly, we use the principal
component analysis (PCA) to get a subspace with dimension $r$.
Then we employ k-means on the PCA-reduced data to get the
clustering results $V'$. Finally, $V$ is initialized as $V=V'+0.3$
and $U$ is by computing the clustering centroid for each class. We
empirically terminate $\ell_1$-NMF, RMF-GMMM, and RMF-LMMM after
$5000$, $500$, and $20$ iterations, respectively.  The clustering
results are shown in the last three rows of Table
\ref{Clustering results}. We can see that RMF-LMMM achieves
tremendous improvements over the two majorant algorithms across
all datasets. RMF-GMMM is also better than $\ell_{1}$-NMF. The lowest relative error in the third row shows that RMF-LMMM can always approximate the measurement matrix much better than the other two.
\renewcommand{\arraystretch}{1.1}
\begin{table}[ht]
\caption{Dataset descriptions, relative errors, and clustering results. For all the
three metrics, larger values are better.}
\label{Clustering results}
\centering
\scriptsize
\begin{tabular}{|c|c|c|c|c|c|}
\hline      
                                               \multicolumn{2}{|c|}{Dataset}             & \multicolumn{1}{c|}{AT\&T} & \multicolumn{1}{c|}{UMINT} & \multicolumn{1}{c|}{CMUPIE} & \multicolumn{1}{c|}{AR} \\ \hline \hline
\multicolumn{1}{|c|}{\multirow{3}{*}{Desc}} & \# Size               & $400$                      & $360$                      & $680$                       & $1300$                  \\
\multicolumn{1}{|c|}{}                             & \# Dim                & $644$                      & $625$                      & $576$                       & $540$                   \\
\multicolumn{1}{|c|}{}                             & \# Class              & $40$                       & $20$                       & $68$                        & $100$                   \\ \cline{1-6} \hline \hline 
\multirow{3}{*}{Error}                             & $L_1$-NMF             & $13.10$                    & $14.52$                    & $15.33$                     & $16.96$                 \\
                                                   & RMF-GMMM              & $13.10$                    & $11.91$                    & $14.64$                     & $17.16$                 \\
                                                   & RMF-LMMM              & $\mathbf{9.69}$            & $\mathbf{9.26}$            & $\mathbf{4.67}$             & $\mathbf{7.53}$         \\ \hline \hline
\multirow{3}{*}{ACC}                               & $L_1$-NMF             & $0.5250$                   & $0.7550$                   & $0.2176$                    & $0.1238$                \\
                                                   & RMF-GMMM              & $0.5325$                   & $0.7889$                   & $0.2691$                    & $0.1085$                \\
                                                   & RMF-LMMM              & $\mathbf{0.7250}$          & $\mathbf{0.8333}$          & $\mathbf{0.4250}$             & $\mathbf{0.2015}$       \\ \cline{1-6}
\multirow{3}{*}{NMI}                               & $L_1$-NMF             & $0.7304$                   & $0.8789$                   & $0.5189$                    & $0.4441$                \\
                                                   & RMF-GMMM              & $0.7499$                   & $0.8744$                   & $0.5433$                    & $0.4289$                \\
                                                   & RMF-LMMM              & $\mathbf{0.8655}$          & $\mathbf{0.9012}$        & $\mathbf{0.6654}$           & $\mathbf{0.4946}$       \\ \cline{1-6}
\multirow{3}{*}{PUR}                               & $L_1$-NMF             & $0.5700$                   & $0.8139$                   & $0.2368$                    & $0.1138$                \\
                                                   & RMF-GMMM              & $0.5725$                   & $0.8222$                   & $0.2897$                    & $0.1354$                \\
                                                   & RMF-LMMM              & $\mathbf{0.7500}$          & $\mathbf{0.8778}$          & $\mathbf{0.4456}$           & $\mathbf{0.2192}$       \\ \cline{1-6}
\end{tabular}
\vspace{-1em}
\end{table}

\section{Conclusions}

In this paper, we propose a weaker condition on surrogates in MM,
which enables better approximation of the objective function. Our
relaxed MM is general enough to include the existing MM methods.
In particular, the non-smooth and non-convex objective function
can be approximated directly, which is never done before. Using
the RMF problems as examples, our locally majorant relaxed MM
beats the state-of-the-art methods with margin, in both solution
quality and convergence speed. We prove that 
the iterates converge to stationary points. To our best knowledge,
this is the first convergence guarantee for variants of RMF
without extra assumptions.

\section{Acknowledgements}
Zhouchen Lin is supported by National Basic Research Program of China (973 Program) (grant no. 2015CB352502), National Natural Science Foundation of China  (NSFC) (grant no. 61272341 and 61231002), and Microsoft Research Asia Collaborative Research Program. Zhenyu Zhao is supported by NSFC (Grant no. 61473302). Hongbin Zha is supported by 973 Program (grant no. 2011CB302202).

{
	\footnotesize
	\bibliographystyle{aaai}
	\bibliography{Relaxed_MM}
}
\newpage

\begin{center}
	\LARGE\textbf{Supplementary Material}
\end{center}

\section{Proofs}
\subsection*{Proof of Proposition 1}
Consider minimizing $g_k(\x)-f(\x)$ in the neighbourhood $\|\x-\x_k\|\leq \epsilon$. It reaches the local minimum $0$ at $\x=\x_k$. By Definition $3$, we have 
\begin{equation}\label{gf} 
   \nabla g_k(\x_k;\dd) \geq \nabla f(\x_k; \dd), \quad  \forall~\x_k+\dd \in \mathcal{C}, \|\dd\|< \epsilon.  
\end{equation} 
Denote ${l}_k(\x)$ as 
\begin{equation}
 {l}_k(\x) =\hat{g}_k(\x)- \gamma_l \|\x-\x_k\|^2_2,
\end{equation}
Similarly, we have
\begin{equation}\label{fl}
   \nabla f(\x_k;\dd) \geq \nabla l_k(\x_k; \dd), \quad  \forall~\x_k+\dd \in \mathcal{C}, \|\dd\|< \epsilon.  
\end{equation} 
  Comparing $g_k(\x)$ with $l_k(\x)$, they are combined with two parts, the common part $\hat{f}_k(\x)$ and the continuously differentiable part $\|\x-\x_k\|_2^2$. And the function $g_k(\x)-l_k(\x) = (\gamma_u+\gamma_l)\|\x-\x_k\|^2$ achieves its global minimum at $\x=\x_k$. Hence the first order optimality condition \cite[Assumption~1.(A3)]{uppernonsmooth} implies 
    \begin{equation}\label{gl}
       \nabla g_k(\x_k;\dd) = \nabla l_k( \x_k;\dd) , \quad  \forall~\x_k+\dd \in \mathcal{C}, \|\dd\|< \epsilon.
\end{equation}
Combining \eqref{gf}, \eqref{fl} and \eqref{gl}, we have
        \begin{equation}
       \nabla g_k(\x_k;\dd) = \nabla f( \x_k;\dd) , \quad  \forall~\x_k+\dd \in \mathcal{C}, \|\dd\|< \epsilon.
\end{equation}
\subsection*{Proof of Theorem 1}

Consider the $\rho$-strongly convex surrogate $g_k(\x)$.  As $\x_{k+1}=\arg \min_{\x \in \mathcal{C}} g_k(\x)$, by the definition of strongly convex function \cite[Lemma B.5]{mairalsurrogate}, we have
     \begin{equation}
       g_k(\x_k)-g_k(\x_{k+1})\geq \frac{\rho}{2} \|\x_k-\x_{k+1}\|_2^2.
\end{equation}
Combining with non-increment of the objective function,
\begin{equation}\label{nonincrement}
f(\x_{k+1}) \leq g_k(\x_{k+1})\leq g_k(\x_k)  = f(\x_k),
\end{equation} 
we have
\begin{equation}\label{sufficientdescreaseoff}
  f(\x_k) - f(\x_{k+1}) \geq \frac{\rho}{2} \|\x_k - \x_{k+1}\|^2.
\end{equation}
Consider the surrogate $g_k(\x)=\dot{g}_k(\x)+\rho/2 \|\x-\x_k\|_2^2$, we have
\begin{equation}
 \begin{split}
 & g_k(\x_{k+1})-f(\x_{k+1}) \\ =& \dot{g}_k(\x_{k+1})+ \frac{\rho}{2} \|\x_{k+1}-\x_k\|_2^2-f(\x_{k+1}) \\ \geq & \frac{\rho}{2} \|\x_{k+1}-\x_k\|_2^2,
 \end{split}
\end{equation}
 where the inequality derived from the locally majorant $\dot{g}_k(\x_{k+1})\geq f(\x_{k+1}) $. Combining with \eqref{nonincrement}, we can also get \eqref{sufficientdescreaseoff}. 
     Thus the sequence $\{f(\x_k)\}$ has sufficient descent. 
     
Summing all the inequalities in \eqref{sufficientdescreaseoff} for $k \geq 1$, we have  
\begin{equation}
  +\infty > f(\x_1)-f(\x_k)|_{k \rightarrow + \infty}\geq \frac{\rho}{2} \sum_{k=1}^{+\infty}   \|\x_k - \x_{k+1}\|^2.
\end{equation}
Then we can infer that
\begin{equation}\label{converge0}
      \lim_{k \rightarrow +\infty } \x_{k+1} - \x_{k}= \mathbf{0}.
\end{equation}
As the sequence $\{\x_k\}$ is bounded, hence has accumulation points. For any accumulation point $\x^*$, there exists a subsequence $\{\x_{k_j}\}$ such that $\lim \limits_{j \to \infty} \x_{k_j} = \x^*$. 

Combining conditions \eqref{surrogate1}, \eqref{surrogate2}, and   \eqref{descending}, we have
\begin{equation}
\begin{split}
    g_{k_j}(\x)& \geq g_{k_j}(\x_{k_j+1})\geq f(\x_{k_j+1}) \\ &\geq f(\x_{k_{j+1}})\geq g_{k_{j+1}}(\x_{k_{j+1}}), \quad \forall \x \in \mathcal{C}.
\end{split}
\end{equation}
Letting $j \rightarrow \infty$ in both sides, we obtain at $\x=\x^*$
  \begin{equation}
   \nabla g(\x^*,\dd) \geq 0, \quad \forall \x^*+\dd  \in \mathcal{C}.
  \end{equation}
Combining with condition \eqref{surrogate3}, we have
  \begin{equation}
   \nabla f(\x^*,\dd) \geq 0, \quad \forall \x^*+\dd \in \mathcal{C}.
  \end{equation}
By Definition $3$, we can conclude that $\x^*$ is a stationary point.     
\subsection*{Proof of Proposition 2}
  Denoting $\Delta \u_i^T$ and $\Delta \v_i^T$ as the $i$-th rows of $\Delta U$ and $\Delta V$, respectively, we can relax
     \begin{equation} \label{eq:in2}
      \begin{split}
      &\|W \odot (\Delta U \Delta V ^T)\|_1 \\ = &
      \begin{Vmatrix} W \odot
      \begin{pmatrix}
    \Delta \u_1^T \Delta \v_1  & \dots & \Delta \u_1^T \Delta \v_n\\
    \vdots & \ddots & \vdots\\
   \Delta \u_m^T \Delta \v_1  & \dots & \Delta \u_m^T \Delta \v_n\\
   \end{pmatrix}
   \end{Vmatrix}_1       \\
     &\leq \frac{\bar{\rho}_u}{2}\|\Delta U\|_F^2 + \frac{\bar{\rho}_v}{2}\|\Delta V\|_F^2, \\
\end{split}
  \end{equation}
 where the inequality is derived from the Cauchy-Schwartz inequality with $\bar{\rho}_u=\# W_{(i,.)}+\epsilon,~\forall i=1,\ldots,m$ and $\bar{\rho}_v=\# W_{(.,j)}+\epsilon,~  \forall i=1,\ldots,n$. The equality holds if and only if $(\Delta U, \Delta V)=\mathbf{(0,0)}$. Then we have 
 \begin{equation}
 \begin{split}
  & \hat{G}_k(\Delta U, \Delta V) + \frac{\bar{\rho}_u}{2}\|\Delta U\|_F^2 + \frac{\bar{\rho}_v}{2}\|\Delta V\|_F^2 
  \\\geq &\hat{G}_k(\Delta U, \Delta V)+ \|W \odot (\Delta U \Delta V ^T)\|_1  
  \\\geq & F_k(\Delta U, \Delta V), 
  \end{split}
 \end{equation}
where the second inequality is derived from the triangular inequality of norms. Similarly we have 
  \begin{equation}
 \begin{split}
 & F_k(\Delta U, \Delta V)
  \\\geq & \hat{G}_k(\Delta U, \Delta V) -\|W \odot (\Delta U \Delta V ^T)\|_1  
  \\\geq &\hat{G}_k(\Delta U, \Delta V)-\frac{\bar{\rho}_u}{2}\|\Delta U\|_F^2 - \frac{\bar{\rho}_v}{2}\|\Delta V\|_F^2   
  \end{split}
 \end{equation}
\subsection*{Proof of Theorem 2}
 Combining $F_k(\Delta U, \Delta V)$ with $G_k(\Delta U, \Delta V)$, we can easily get $F_k(\mathbf{(0,0)})=G_k(\mathbf{(0,0)})$. Let $(\Delta U_k, \Delta V_k)$ represent the minimizer of $G_k(\Delta U, \Delta V)$. By the choice of $(\rho_u,\rho_v)$, both RMF-GMMM and RMF-LMMM can ensure 
   \begin{equation}
     G_k(\Delta U_k, \Delta V_k) \geq F_k(\Delta U_k, \Delta V_k).
\end{equation}
Combining Proposition $1$ and $2$, we can ensure the first order smoothness in the infinity
     \begin{equation}
     \begin{split} 
 \lim_{k \to \infty}&\left( \nabla F_k( \mathbf{0, 0};D_u,D_v)-\nabla G_k( \mathbf{0, 0};D_u,D_v) \right)= 0, \\& \quad \forall~ U_k+D_u \in \mathcal{C}_u, V_k+D_v \in \mathcal{C}_v.
 \end{split} 
\end{equation}
In addition, $G_k(\Delta,\Delta V)$ is strongly convex and the sequence $\{(U_k,V_k)\}$ is bounded by the constraints and regularizations, 
according to Theorem $1$, the original function sequence $\{F(U_k,V_k)\}$ would have sufficient descent, and any limit point of the sequence $\{(U_k,V_k)\}$ is a stationary point of the objective function $F(U,V)$.
      
\section{Minimizing the surrogate by LADMPSAP}
\subsection{Sketch of LADMPSAP}
LADMPSAP fits for solving the following linearly constrained separable
convex programs:
\begin{equation}
\min \limits_{\x_1,\cdots,\x_n} \sum\limits_{j=1}^n f_j(\x_j),\quad
s.t.\quad
\sum \limits_{j=1}^n \mathcal{A}_j(\x_j)=\mathbf{b},\label{eq:separable constraint2}
\end{equation}
where $\x_j$ and $\mathbf{b}$ could be either vectors or
matrices, $f_j$ is a proper convex function, and $\mathcal{A}_j$ is a linear mapping. Very often, there are multiple blocks of variables ($n\geq 3$).
We denote the iteration index by superscript $i$. The LADMPSAP algorithm consists of the following steps \citep{ladmpsap}:
\begin{enumerate}[(a)]
\item Update $\x_j$'s ($j=1,\cdots,n$) in parallel:\\
\begin{equation}
\begin{split}
\x_j^{i+1}=\argmin\limits_{\x_j} f_j(\x_j)
+ &\dfrac{\sigma_j^{(i)}}{2} \lbar\z_j-\x_j^{i}+\A_j^{\dag}(\hat{\y}^{i})/\sigma_j^{(i)}\rbar^2 \label{eq:update_xi_ladmpsap}
\end{split}
\end{equation}
\item Update $\y$:\\
\begin{equation}
\y^{i+1} = \y^i +
\beta^{(i)}\left(\sum\limits_{j=1}^n\A_j(\x_j^{i+1})-\b\right).
\label{eq:update_y_ladmpsap}
\end{equation}
\item Update $\beta$:\\
\begin{equation}
\beta^{(i+1)} = \min(\beta^{\max},\rho\beta^{(i)}),
\label{eq:update_beta_ladmpsap}
\end{equation}
\end{enumerate}
 where $\y$ is the Lagrange multiplier, $\beta^{(i)}$ is the penalty parameter, $\beta^{max}\gg 1$ is an upper bound of $\beta^{(i)}$, $\sigma_j^{(i)}=\eta_j\beta^{(i)}$ with $\eta_j
> n\|\A_j\|^2$ ($\|\A_j\|$ is the operator norm of $\A_j$), $\A_j^{\dag}$ is the adjoint operator of $\A_j$,
\begin{equation}
\hat{\y}^{i}=\y^i +
\beta^{(i)}\left(\sum\limits_{j=1}^n\A_j(\x_j^{i})-\b\right),\label{eq:hat_lambda}
\end{equation}
and
\begin{equation}
\rho = \left\{\begin{array}{ll} \rho_0, &  \mbox{if} \
\beta^{(i)}\max\left(\left\{\sqrt{\eta_j}\lbar\z_j^{i+1}-\x_j^i\rbar\right\}  \right)  /\lbar\mathbf{b}\rbar < \varepsilon_1,
 \\
1,  & \mbox{otherwise},
\end{array}\right.\label{eq:update_rho_ladmpsap}
\end{equation}
with $\rho_0 \geq 1$ being a constant and $0 < \varepsilon_1\ll 1$
being a threshold. 
The iteration terminates when the following two conditions are met:
\begin{equation}
\beta_{k}\max\left(\left\{\sqrt{\eta_i}\lbar\mathbf{x}_i^{k+1}-\mathbf{x}_i^k\rbar,i=1,\cdots,n\right\}\right)/\|\mathbf{b}\|
< \varepsilon_1,\label{eq:stopping1}
\end{equation}
\begin{equation}
\left\|\sum\limits_{i=1}^n\mathcal{A}_i(\mathbf{x}_i^{k+1}) -
\mathbf{b}\right\|/\|\mathbf{b}\| <
\varepsilon_2.\label{eq:stopping2}
\end{equation}


\subsection{The Optimization Using LADMPSAP}
We aim to minimize 
  \begin{equation}\label{eq:eqconstraint2}
   \begin{split}
 &\min_{E,\Delta U,\Delta V}  \|W \odot E\|_1 \\ 
 &+\left(\frac{\rho_u}{2}\|\Delta U\|_F^2+R_u(U_k+\Delta U)+ \delta_{\mathcal{C}_u}(U_k+\Delta U)\right)
 \\& + \left(\frac{\rho_v}{2}\|\Delta V\|_F^2+R_v(V_k+\Delta V)+ \delta_{\mathcal{C}_v}(V_k+\Delta V)\right), \\
                  &\text{s.t.  } E+\Delta U V_k^T+U_k\Delta V^T=M-U_kY_k^T,
              \end{split}
              \end{equation}
 where the function $\delta_\mathcal{C}(\x): \mathbb{R}^p \rightarrow \mathbb{R}$ is define as:
  \begin{equation}
     \delta_{\mathcal{C}}(\x)=\left\{ \begin{array}{ll}
     0,  &\quad  \mbox{if } \x \in \mathcal{C},  \\
     +\infty, &\quad \mbox{otherwise.}
     \end{array} 
        \right.
\end{equation}     
  which naturally fits into the model problem \eqref{eq:separable constraint2}.

  According to \eqref{eq:update_xi_ladmpsap}, $E$ can be updated by solving the following subproblem
  \begin{equation}\label{eq:eproblem}
  \min_E  \|W \odot E\|_1 +\frac{\sigma_e^{(i)} }{2}\|E-E^i+\hat{Y}^i/\sigma_e^{(i)}\|_F^2,
\end{equation}
 where
    \begin{equation}
 \hat{Y}^i=Y^i+\beta^{(i)}( E^i + \Delta U^i V_k^T+U_k \Delta V^{iT}+U_kY_k^T-M),
 \end{equation}
 and $\sigma_e^{i}=\eta_e\beta^{(i)}$. We choose $\eta_e=3L_e+\epsilon$, where 3 is the number of parallelly updating variables, i.e., $E$, $\Delta U$ and $\Delta V$, $L_e$ denotes the squared spectral norm of the linear mapping on $E$, which is equal to $1$, and $\epsilon$ is a small positive scalar. The solution to \eqref{eq:eproblem} is
  \begin{equation}\label{eq:updatee}
   \begin{split}
   E^{i+1} = & W \odot \mathcal{S}_{{\sigma_e^{(i)}}^{-1}} \bigl( E^i-\hat{Y}^i/\sigma_e^{(i)}\bigr)
    \\& + \bar{W}\odot \bigl(E^i-\hat{Y}^i/\sigma_e^{(i)}\bigr),
   \end{split}
  \end{equation}
 where $\mathcal{S}$ is the shrinkage operator \citep{zlinalm}:
     \begin{equation}
       \mathcal{S}_\gamma(x)=\max(|x|-\gamma,0)\mbox{sgn}(x),
     \end{equation}
 and $\bar{W}$ is the complement of $W$.

 Also by \eqref{eq:update_xi_ladmpsap}, $\Delta U$ and $\Delta V$ can be updated by solving:
   \begin{equation}\label{uproblem}
     \begin{split}
     \min_{\Delta U}&
     \frac{\rho_u}{2}\|\Delta U\|_F^2+R_u(U_k+\Delta U)+ \delta_{\mathcal{C}_u}(U_k+\Delta U)  
    \\& +\frac{\sigma_u^{(i)}}{2}\|\Delta U -\Delta U^i+  \hat{Y}^i V_k/\sigma_u^{(i)}  \|_F^2,
    \end{split}
   \end{equation}
      \begin{equation}\label{vproblem}
     \begin{split}
     \min_{\Delta V}&
     \frac{\rho_v}{2}\|\Delta V\|_F^2+R_v(V_k+\Delta V)+ \delta_{\mathcal{C}_v}(V_k+\Delta V) 
     \\& +\frac{\sigma_v^{(i)}}{2}\|\Delta V -\Delta V^i+ \hat{Y}^{iT} U_k /(\sigma_v^{(i)})\|_F^2,
    \end{split}
   \end{equation}
   where $\sigma_u^{(i)}=\eta_x \beta^{(i)}$, $\eta_x=3\|V_k\|^2_2+\epsilon$,  $\sigma_v=\eta_v \beta^{(i)}$, and $\eta_v=3\|U_k\|^2_2+\epsilon$ ($\lbar\cdot\rbar_2$ denotes the largest singular value of a matrix). 
 
 For \textbf{Low Rank Matrix Recovery}, \eqref{uproblem} and \eqref{vproblem} can be solved by simply taking the derivative w.r.t $\Delta U$ and $\Delta V$ to be $\mathbf{0}$,
\begin{equation}
       \begin{split}
    \Delta U^{i+1} = &\left(-\lambda_u U_k+ \sigma_u^{(i)} \Delta U^i-
     \hat{Y}^i V_k\right)/(\lambda_u +\sigma_u ^{(i)} +\rho_u ),
    \end{split}
   \end{equation}
   \begin{equation}
       \begin{split}
    \Delta V^{i+1} =  & \left(-\lambda_v V_k+ \sigma_v^{(i)} \Delta V^i-\hat{Y}^{iT} U_k \right)/(\lambda_v +\sigma_v^{(i)}  +\rho_v).
    \end{split}
   \end{equation}  

For \textbf{Non-negative Matrix Factorization},  \eqref{uproblem} and \eqref{vproblem} can be solved by projecting on to the non-negative subspace,
\begin{equation}
\begin{aligned}
   \Delta U^{i+1}= \mathcal{S}^+_0
    &\bigl((\sigma_u ^{(i)} + \rho_u) U_k+ \sigma_u^{(i)} \Delta U^i \\&-
     \hat{Y}^i V_k \bigr)/(\lambda_u +\sigma_u ^{(i)} +\rho_u )-U_k,
\end{aligned}
\end{equation}
\begin{equation}
\begin{aligned}
   \Delta V^{i+1}= \mathcal{S}^+_{\lambda_v/(\rho_v+\sigma_v^{(i)})} &\bigl(V_k-\sigma_v^{(i)}/(\rho_v+\sigma_v^{(i)})(-\Delta V^i \\&+ \hat{Y}^{iT} U_k/\sigma_v^{(i)}) \bigr)-V_k,
\end{aligned}
\end{equation}
 where $\mathcal{S}^+$ is the positive shrinkage operator:
     \begin{equation}
       \mathcal{S}_\gamma(x)=\max(x-\gamma,0).
     \end{equation}
     Next, we update $Y$ as \eqref{eq:update_y_ladmpsap}:
\begin{equation}\label{eq:updatey}
  \begin{split}
   Y^{i+1}=Y^i+\beta^{(i)}&(E^{i+1}  +\Delta U^{i+1} V_k^T \\& +U_k \Delta V^{(i+1)T}  +U_kV_k^T-M),
  \end{split}
\end{equation}
and update $\beta$ as \eqref{eq:update_beta_ladmpsap}:
\begin{equation} \label{beta1}
    \beta^{(i+1)}=\min(\beta^{\max},\rho \beta^{(i)}),
\end{equation}
 where $\rho$ is defined as \eqref{eq:update_rho_ladmpsap}:
\begin{equation} \label{beta2}
 \rho = \left\{\begin{array}{l}
  \begin{aligned}   \rho_0, & \quad  \textrm{if } \beta^{(i)} \max \bigl(\sqrt{\eta_e}\|E^{i+1}-E^{i}\|_F, \\& \sqrt{\eta_u}\|\Delta U^{i+1}-\Delta U^{i}\|_F, \\ & \sqrt{\eta_v}\|\Delta V^{i+1}-\Delta V^{i}\|_F \bigr)/\|M-U_kV_k^T\|_F < \varepsilon_1,
   \end{aligned} \\
 1, \quad \text{ otherwise.}
 \end{array}\right.
\end{equation}
We terminate the iteration when the following two conditions are met:
   \begin{equation}\label{stop1}
   \begin{split}
   \beta^{(i)} \max(&\sqrt{\eta_e}\|E^{i+1}-E^{i}\|_F, \sqrt{\eta_u}\|\Delta U^{i+1}-\Delta U^{i}\|_F, \\ & \sqrt{\eta_v}\|\Delta V^{i+1}-\Delta V^{i}\|_F)/\|M-U_kV_k^T\|_F < \varepsilon_1
   \end{split}
   \end{equation}
   \begin{equation}\label{stop2}
     \|E^{i+1}-\Delta U^{i+1} V_k^T-U_k\Delta V^{(i+1)T}\|_F/\|M-U_kV_k^T\|_F < \varepsilon_2.
   \end{equation}
   For better reference, we summarize the algorithm for minimizing $G_k(\Delta U, \Delta V)$ in Algorithm \ref{alg:ladmpsap}. When first executing Algorithm \ref{alg:ladmpsap}, we initialize $E^0=M-U_0V_0^T$, $\Delta U^0=\mathbf{0}$, $\Delta V^0= \mathbf{0}$ and $Y^0=\mathbf{0}$. In the subsequent main iterations, we adopt the warm start strategy. Namely, we initialize $E^0$, $\Delta U^0,\Delta V^0$ and $Y^0$ with their respective optimal values in last main iteration.
\subsection{Parameter Setting for LADMPSAP}
\textbf{Low Rank Matrix Recovery}: we set
$\varepsilon_1=10^{-5}$, $\varepsilon_2=10^{-4}$, $\rho_0=1.5$,
and $\beta^{\max}=10^{10}$ as the default value.

\vspace{0.4em}  
\noindent \textbf{Non-negative Matrix Factorization}: we set $\varepsilon_1=10^{-4}$, $\varepsilon_2=10^{-4}$, $\rho_0=3$,
and $\beta^{\max}=10^{10}$ as the default value.
    \renewcommand{\algorithmicrequire}{\textbf{Input:}}
 \renewcommand{\algorithmicensure}{\textbf{Output:}}
\begin{algorithm}[htpb]
 \caption{Minimizing $G_k(\Delta U, \Delta V)$  via LADMPSAP}
 \label{alg:ladmpsap}
\begin{algorithmic}[1]
   \STATE  Initialize $i=0$, $E^0$, $\Delta U^0$, $\Delta V^0$, $Y^0$, $\rho_0>1$ and $\beta^0\propto (m+n)\varepsilon_1$.
\WHILE{ \eqref{stop1} or \eqref{stop2} is not satisfied}
\STATE Update $E$, $\Delta U$, and $\Delta V$ parallelly, accordingly. 
\STATE Update $Y$ as \eqref{eq:updatey}.
\STATE Update $\beta$ as \eqref{beta1} and \eqref{beta2}.
\STATE $i=i+1$.
\ENDWHILE
\ENSURE The optimal $(\Delta U_k, \Delta V_k)$.
\end{algorithmic}
\end{algorithm}


\end{document}